\documentclass[reqno,11pt]{amsart}
\usepackage{amscd,amssymb}

\numberwithin{equation}{section}

\theoremstyle{plain}
\newtheorem{Thm}[subsection]{Theorem}
\newtheorem{Cor}[subsection]{Corollary}
\newtheorem{Lem}[subsection]{Lemma}
\newtheorem{Prop}[subsection]{Proposition}

\theoremstyle{definition}
\newtheorem{Def}[subsection]{Definition}

\theoremstyle{remark}

\newtheorem{Rem}[subsection]{Remark}

\newif\ifShowLabels
\ShowLabelstrue
\newdimen\theight
\def\TeXref#1{%
        \leavevmode\vadjust{\setbox0=\hbox{{\tt
                \quad\quad  {\small \textrm #1}}}%
        \theight=\ht0
        \advance\theight by \lineskip
        \kern -\theight \vbox to
        \theight{\rightline{\rlap{\box0}}%
        \vss}%
        }}%

\ShowLabelsfalse

\renewcommand{\sec}[2]{\section{#2}\label{S:#1}%
        \ifShowLabels \TeXref{{S:#1}} \fi}
\newcommand{\ssec}[2]{\subsection{#2}\label{SS:#1}%
        \ifShowLabels \TeXref{{SS:#1}} \fi}

\newcommand{\refs}[1]{Section ~\ref{S:#1}}
\newcommand{\refss}[1]{Subsection ~\ref{SS:#1}}
\newcommand{\reft}[1]{Theorem ~\ref{T:#1}}
\newcommand{\refl}[1]{Lemma ~\ref{L:#1}}
\newcommand{\refp}[1]{Proposition ~\ref{P:#1}}
\newcommand{\refc}[1]{Corollary ~\ref{C:#1}}

\newcommand{\refr}[1]{Remark ~\ref{R:#1}}
\newcommand{\refe}[1]{\eqref{E:#1}}

\newenvironment{thm}[1]%
        { \begin{Thm} \label{T:#1}  \ifShowLabels \TeXref{T:#1} \fi }%
        { \end{Thm} }

\renewcommand{\th}[1]{\begin{thm}{#1} \sl }
\renewcommand{\eth}{\end{thm} }

\newenvironment{lemma}[1]%
        { \begin{Lem} \label{L:#1}  \ifShowLabels \TeXref{L:#1} \fi }%
        { \end{Lem} }
\newcommand{\lem}[1]{\begin{lemma}{#1} \sl}
\newcommand{\elem}{\end{lemma}}

\newenvironment{propos}[1]%
        { \begin{Prop} \label{P:#1}  \ifShowLabels \TeXref{P:#1} \fi }%
        { \end{Prop} }
\newcommand{\prop}[1]{\begin{propos}{#1}\sl }
\newcommand{\eprop}{\end{propos}}

\newenvironment{corol}[1]%
        { \begin{Cor} \label{C:#1}  \ifShowLabels \TeXref{C:#1} \fi }%
        { \end{Cor} }
\newcommand{\cor}[1]{\begin{corol}{#1} \sl }
\newcommand{\ecor}{\end{corol}}

\newenvironment{defeni}[1]%
        { \begin{Def} \label{D:#1}  \ifShowLabels \TeXref{D:#1} \fi }%
        { \end{Def} }
\newcommand{\defe}[1]{\begin{defeni}{#1} \sl }
\newcommand{\edefe}{\end{defeni}}

\newenvironment{remark}[1]%
        { \begin{Rem} \label{R:#1}  \ifShowLabels \TeXref{R:#1} \fi }%
        { \end{Rem} }
\newcommand{\rem}[1]{\begin{remark}{#1}}
\newcommand{\erem}{\end{remark}}

\newcommand{\eq}[1]%
        { \ifShowLabels \TeXref{E:#1} \fi
           \begin{equation} \label{E:#1} }
\newcommand{\eeq}{\end{equation}}

\newcommand{\prf}{ \begin{proof} }
\newcommand{\eprf}{ \end{proof} }

\newcommand\alp{\alpha}         
\newcommand\bet{\beta}
\newcommand\gam{\gamma}         \newcommand\Gam{\Gamma}
\newcommand\del{\delta}         \newcommand\Del{\Delta}
\newcommand\eps{\varepsilon}

\newcommand\iot{\iota}

\newcommand\lam{\lambda}                \newcommand\Lam{\Lambda}
\newcommand\sig{\sigma}

\newcommand\calA{{\mathcal{A}}}

\newcommand\calC{{\mathcal{C}}}

\newcommand\calE{{\mathcal{E}}}
\newcommand\calF{{\mathcal{F}}}

\newcommand\calL{{\mathcal{L}}}

\newcommand\calO{{\mathcal{O}}}

\newcommand\calS{{\mathcal{S}}}

\newcommand\calW{{\mathcal{W}}}

\newcommand\calZ{{\mathcal{Z}}}

\newcommand\bfc{{\mathbf c}}

\newcommand\ZZ{\mathbb{Z}}

\newcommand\CC{\mathbb{C}}

 \newcommand\grg{{\mathfrak{g}}}
 \newcommand\grh{{\mathfrak{h}}}

 \newcommand\grk{{\mathfrak{k}}}

\newcommand\nek{,\ldots,}
\newcommand\sdp{\times \hskip -0.3em {\raise 0.3ex
\hbox{$\scriptscriptstyle |$}}} 

\newcommand\Dom{\operatorname{Dom}}

\newcommand\End{\operatorname{End\,}}

\newcommand\Hom{\operatorname {Hom}}

\newcommand\IM{\operatorname{Im}}

\newcommand\Ker{\operatorname{Ker}}

\newcommand\RE{\operatorname{Re}}

\newcommand\spin{\operatorname{spin}}

\newcommand\supp{\operatorname{supp}}

\newcommand\Tr{\operatorname{Tr}}

\newcommand\oB{{\overline{B}}}

\renewcommand\oe{{\overline{e}}}

\newcommand\oz{{\overline{z}}}

\newcommand\tilD{{\widetilde{D}}}

\newcommand\tilJ{{\widetilde{J}}}

\newcommand\tilDel{{\widetilde{\Delta}}}

\newcommand\tileps{{\widetilde{\eps}}}

\renewcommand{\>}{\rangle}
\newcommand{\<}{\langle}

\renewcommand{\d}{\text{\( \partial\)}}
\newcommand{\p}{\bar{\d}}

\newcommand{\n}{\nabla}
\newcommand{\tD}{\tilDel}
\newcommand{\tJ}{\tilJ}

\newcommand{\nk}{\nabla^{\E\otimes\L^k}}
\newcommand{\tilnk}{\widetilde{\nabla}^{\E\otimes\L^k}}
\newcommand{\tildk}{\widetilde{\Del}_k}
\newcommand{\Ek}{{\E\otimes\L^k}}
\newcommand{\E}{\calE}\newcommand{\W}{\calW}
\newcommand{\Z}{\calZ}\newcommand{\F}{\calF}
\newcommand{\A}{\calA}\renewcommand{\O}{\calO}
\renewcommand{\L}{\calL}

\newcommand{\Lk}{\L^k}
\newcommand{\fes}{F^{\E/\calS}}

\newcommand{\EndC}{\operatorname{End}_{C(M)}\,}
\newcommand{\g}{{\Gam}}
\newcommand{\gc}{{\Gam(M,C(M))}}
\newcommand{\gme}{{\Gam(M,\E)}}
\newcommand{\ha}{^{1,0}}
\newcommand{\ah}{^{0,1}}

\newcommand{\hor}{^{\text{hor}}}
\renewcommand{\vert}{^{\text{vert}}}

\newcommand{\G}{\Gamma}

\newcommand{\ka}{K\"ahler }
\begin{document}

\title{Vanishing theorems on covering manifolds}
\author{Maxim Braverman}
\address{Department of Mathematics\\
        The Ohio State University   \\
        Columbus, OH 43210 \\
        USA
         }
\email{maxim@math.ohio-state.edu}
\thanks{This research was partially supported by grant No. 96-00210/1 from
the United States-Israel Binational Science Foundation (BSF)}

\begin{abstract}
Let $M$ be an oriented even-dimensional Riemannian manifold on which a discrete
group $\G$ of orientation-preserving isometries acts freely, so that the
quotient $X=M/\G$ is compact. We prove a vanishing theorem for a half-kernel of
a $\G$-invariant Dirac operator on a $\G$-equivariant Clifford module over $M$,
twisted by a sufficiently large power of a $\G$-equivariant line bundle, whose
curvature is non-degenerate at any point of $M$. This generalizes our previous
vanishing theorems for Dirac operators on a compact manifold.

In particular, if $M$ is an almost complex manifold we prove a vanishing theorem
for the half-kernel of a $\spin^c$ Dirac operator, twisted by a line bundle with
curvature of a mixed sign. In this case we also relax the assumption of
non-degeneracy of the curvature. When $M$ is a complex manifold our results
imply analogues of Kodaira and Andreotti-Grauert vanishing theorems for covering
manifolds.

As another application, we show that semiclassically the $\spin^c$ quantization
of an almost complex covering manifold gives an ``honest" Hilbert space. This
generalizes a result of Borthwick and Uribe, who considered quantization of
compact manifolds.

Application of our results to homogeneous manifolds of a real semisimple Lie
group leads to new proofs of Griffiths-Schmidt and Atiyah-Schmidt vanishing
theorems.
\end{abstract}
\maketitle
\tableofcontents
\sec{introd}{Introduction}

One of the most fundamental results of the geometry of {\em compact} complex
manifolds is the Kodaira vanishing theorem for the cohomology of the sheaf of
sections of a holomorphic vector bundle twisted by a large power of a positive
line bundle. Andreotti and Grauert \cite{AndGr62} generalized this result to the
case when the line bundle is not necessarily positive (but satisfies some
non-degeneracy conditions, cf. \refs{cvanish}).

Both the Kodaira and the Andreotti-Grauert theorems are equivalent to a
vanishing of the kernel of the restriction of the Dolbeault-Dirac operator
$\sqrt2(\p+\p^*)$ to the space of differential forms of certain degree. In
\cite{BrVanish}, the author obtained a generalization of these results to
abstract Dirac operators twisted by a large power of a line bundle (see also
\cite{BoUr96} where an analogue of the Kodaira vanishing theorem for $\spin^c$
Dirac operator on an almost \ka manifold is proven).

In this paper we show that suitable generalizations of all the preceding
vanishing theorems remain true if the base manifold is not compact but is an
infinite normal covering of a compact manifold. The obtained results are very
convenient for the study of homogeneous vector bundles over homogeneous spaces
of real semisimple Lie groups. In particular, we obtain new proofs of certain
results of Griffiths and Schmid \cite{GrifSch69} and Atiyah and Schmid
\cite{AtSch77}, cf. \refs{repr}.

As another application, we prove that {\em semiclassically the
$\spin^c$-quantization of an almost complex covering manifold gives an honest
Hilbert space} and not just a virtual one, cf. \refss{BU}.

We now give a brief review of the main results of the paper.

\ssec{DV}{The $L^2$ vanishing theorem for the half-kernel of a Dirac operator}
Suppose $M$ is an oriented even-dimensional Riemannian manifold on which a
discrete group $\G$ of orientation-preserving isometries acts freely, so that
the quotient $M/\G$ is compact. Let $C(M)$ denote the Clifford bundle of $M$,
i.e., a vector bundle whose fiber at any point is isomorphic to the Clifford
algebra of the cotangent space. Let $\E$ be a $\G$-equivariant self-adjoint
Clifford module over $M$, i.e., a $\G$-equivariant vector bundle over $M$
endowed with a $\G$-invariant Hermitian structure and a fiberwise self-adjoint
action of $C(M)$. Then (cf. \refss{chir}) $\E$ possesses a natural grading
$\E=\E^+\oplus\E^-$. Let $\L$ be a $\G$-equivariant Hermitian line bundle
endowed with a $\G$-invariant Hermitian connection $\n^\L$. These data define
(cf. \refss{dirac}) a $\G$-invariant self-adjoint Dirac operator $D_k$ acting on
the space $L^2(M,\Ek)$ of square integrable sections of $\Ek$. The curvature
$F^\L$ of $\n^\L$ is an imaginary valued 2-form on $M$. If it is non-degenerate
at all points of $M$, then $iF^\L$ is a symplectic form on $M$, and, hence,
defines an orientation of $M$. Our first result (\reft{Dirvanish}) states that
{\em the restriction of the kernel of $D_k$ to $L^2(M,\E^-\otimes\Lk)$ (resp. to
$L^2(M,\E^+\otimes\Lk)$) vanishes for large $k$ if this orientation coincides
with (resp. is opposite to) the given orientation of $M$.}

\ssec{ag}{The $L^2$ Andreotti-Grauert theorem}
Suppose that a discrete group $\G$ acts holomorphically and freely of a complex
manifold $M$, so that the quotient $M/\G$ is compact. Let $\W$ be a holomorphic
$\G$-equivariant vector bundle over $M$ and let $\L$ be a holomorphic
$\G$-equivariant line bundle over $M$. Assume that $\L$ carries a $\G$-invariant
Hermitian metric whose curvature form has at least $q$ negative and at least $p$
positive eigenvalues at any point $x\in M$. Then (\reft{AG}), {\em the
$L^2$-cohomology $L^2H^{0,j}(M,\W\otimes\Lk)$ of $M$ with coefficients in the
tensor product $\W\otimes\Lk$ vanishes for $j\not=q, q+1\nek n-p$ and $k\gg0$.}

In particular, {\em if $\L$ is a positive bundle, then
$L^2H^{0,j}(M,\W\otimes\Lk)=0$, for $j\not=0$ and $k\gg0$}. This is an
$L^2$-analogue of the Kodaira vanishing theorem.

\ssec{ac}{A generalization to almost complex manifolds}
If, in the conditions of the previous subsection, $M$ is a \ka manifold, then
the $L^2$ cohomology $L^2H^{0,*}(M,\W\otimes\Lk)$ is isomorphic to the kernel of
the Dolbeault-Dirac operator
$D_k:L^2\A^{0,*}(M,\W\otimes\Lk)\to{}L^2\A^{0,*}(M,\W\otimes\Lk)$, where
$L^2\A^{0,*}(M,\W\otimes\Lk)$ denotes the space of square integrable
differential forms of type $(0,*)$ on $M$ with coefficients in $\W\otimes\Lk$.
This suggest a generalization of the $L^2$ Andreotti-Grauert theorem to the case
when $M$ is only an almost complex manifold.

Assume that $\L$ possess a Hermitian connection whose curvature is a $(1,1)$
form on $M$ which has at least $q$ negative and at least $p$ positive
eigenvalues at any point $x\in M$. In this situation a Dirac operator
$D_k:L^2\A^{0,*}(M,\W\otimes\Lk)\to{}L^2\A^{0,*}(M,\W\otimes\Lk)$ is defined,
cf. \refss{acAG}. If the almost complex structure on $M$ is not integrable, one
can not hope that the kernel of $D_k$ belongs to
$\oplus_{j=q}^{n-p}L^2\A^{0,j}(M,\W\otimes\Lk)$. However, we show in \reft{UB+},
that {\em for any $k\gg0$ and any $\alp\in\Ker D_k$, ``most of the norm" of
$\alp$ is concentrated in $\oplus_{j=q}^{n-p}L^2\A^{0,j}(M,\W\otimes\Lk)$}. In
particular, {\em if the curvature of $\L$ is non-degenerate and has exactly $q$
negative eigenvalues at any point of $M$, then ``most of the norm" of
$\alp\in\Ker D_k$ is concentrated in $L^2\A^{0,q}(M,\W\otimes\Lk)$, and,
depending on the parity of $q$, the restriction of the kernel of $D_k$ either to
$L^2\A^{0,\text{odd}}(M,\W\otimes\Lk)$ or to
$L^2\A^{0,\text{even}}(M,\W\otimes\Lk)$ vanishes.}

In \refss{BU}, we discuss applications of the above result to geometric
quantization of covering manifolds. In this way we obtain $L^2$ analogues of
results of Borthwick and Uribe \cite{BoUr96}.

\subsection*{Contents}
The paper is organized as follows:

In Sections \ref{S:Dirvanish}--\ref{S:almcomp}, we state our vanishing theorems
for covering manifold.

In \refs{repr}, we discuss applications of our results to representation theory
of real semisimple Lie groups.

In \refs{prdirvanish}, we present the proof of \reft{Dirvanish} (the vanishing
theorem for the kernel of a Dirac operator). The proof is based on two
statements (Propositions~\ref{P:D-D} and \ref{P:lapl}) which are proven in the
later sections.

In \refs{prUB}, we prove an estimate on the Dirac operator on an almost complex
manifold (\refp{estDk}) and use it to prove \reft{UB+} (our analogue of the
Andreotti-Grauert vanishing theorem for almost complex manifolds). The proof is
based on Propositions~\ref{P:lapl} and \ref{P:D-Dac} which are proven in later
sections.

In \refs{prAG}, we prove the $L^2$ Andreotti-Grauert theorem (\reft{AG}).

In \refs{compare}, we use the Lichnerowicz formula to prove
Propositions~\ref{P:D-D}, \ref{P:D-Dac} and \ref{P:gr}. These results establish
the connection between the Dirac operator and the rough Laplacian. They are used
in the proofs of Theorems~\ref{T:Dirvanish}, \ref{T:UB+} and \ref{T:AG}.

{}Finally, in \refs{lapl}, we prove \refp{lapl} (the estimate on the rough
Laplacian).

\subsection*{Acknowledgments}
I would like to thank Yael Karshon for valuable discussions.

\sec{Dirvanish}{$L^2$ vanishing theorem for the half-kernel of a Dirac operator}

In this section we formulate one of the main results of the paper: the $L^2$
vanishing theorem for the half-kernel of a Dirac operator (cf.
\reft{Dirvanish}).

The section is organized as follows: in
Subsections~\ref{SS:clmodule}--\ref{SS:dirac} we recall some basic facts about
Clifford modules and Dirac operators. When possible we follow the notations of
\cite{BeGeVe}. In \refss{Gamma} we discuss some properties of Dirac operator on
covering manifolds. Finally, in \refss{Dirvanish} we formulate our main result.

\ssec{clmodule}{Clifford Modules}
Suppose $M$ is an oriented even-dimensional Riemannian manifold and let $C(M)$
denote the Clifford bundle of $M$ (cf. \cite[\S3.3]{BeGeVe}), i.e., a vector
bundle whose fiber at any point $x\in M$ is isomorphic to the Clifford algebra
$C(T^*_xM)$ of the cotangent space.

A {\em Clifford module \/} on $M$ is a complex vector bundle $\E$ on $M$ endowed
with an action of the bundle $C(M)$. We write this action as
\[
        (a,s) \ \mapsto \ c(a)s, \quad \mbox{where} \quad
                        a\in \gc, \ s\in \gme.
\]

A Clifford module $\E$ is called {\em self-adjoint \/} if it is endowed with a
Hermitian metric such that the operator $c(v):\E_x\to\E_x$ is skew-adjoint, for
any $x\in M$ and any $v\in T_x^*M$.

A connection $\n^\E$ on a Clifford module $\E$ is called a {\em Clifford
connection \/} if
\[
        [\n^\E_X,c(a)] \ = \ c(\n_X a), \quad
                \mbox{for any} \quad  a\in \gc, \ X\in\g(M,TM).
\]
In this formula, $\n_X$ is the Levi-Civita covariant derivative on $C(M)$
associated with the Riemannian metric on $M$.

Suppose $\E$ is a Clifford module and $\calW$ is a vector bundle over $M$. The
{\em twisted Clifford module obtained from $\E$ by twisting with $\calW$ \/} is
the bundle $\E\otimes\calW$ with Clifford action $c(a)\otimes1$. Note that the
twisted Clifford module $\E\otimes\calW$ is self-adjoint if and only if so is
$\E$.

Let $\n^\calW$ be a connection on $\calW$ and let $\n^\E$ be a Clifford
connection on $\E$. Then the {\em product connection}
\begin{equation}\label{E:nEW}
        \n^{\E\otimes\calW} \ = \ \n^\E\otimes 1 \ + \ 1\otimes \n^\calW
\end{equation}
is a Clifford connection on $\E\otimes\calW$.

\ssec{chir}{The chirality operator. The natural grading}
Fix $x\in M$ and let $e_1\nek e_{2n}$ be an oriented orthonormal basis of
$T_x^*M$. Consider the element
\begin{equation}\label{E:Gam}
        \gam \ = \ i^n\, e_1\cdots e_{2n} \ \in \ C(T_x^*M)\otimes\CC,
\end{equation}
called the {\em chirality operator}. It is independent of the choice of the
basis, anti-commutes with any $v\in T_x^*M\subset C(T_x^*M)$, and satisfies
$\gam^2=-1$, cf. \cite[\S3.2]{BeGeVe}. We also denote by $\gam$ the section of
$C(M)$ whose restriction to each fiber is equal to the chirality operator.

The {\em natural grading \/} on a Clifford module $\E$ is defined by the formula
\begin{equation}\label{E:grad}
        \E^\pm \ = \ \{v\in \E: \ c(\gam)\, v=\pm v\}.
\end{equation}
Note that this grading is preserved by any Clifford connection on $\E$. Also, if
$\E$ is a self-adjoint Clifford module (cf. \refss{clmodule}), then the
chirality operator $c(\gam):\E\to\E$ is self-adjoint. Hence, the subbundles
$\E^\pm$ are orthogonal with respect to the Hermitian metric on $\E$. {\em In
this paper we endow all our Clifford modules with the natural grading.}

\ssec{dirac}{Dirac operators}
The {\em Dirac operator \/} $D:\gme\to\gme$ associated to a Clifford connection
$\n^\E$ is defined by the following composition
\begin{equation}\label{E:dir1}
\begin{CD}
        \gme @>\n^\E>> \g(M,T^*M\otimes \E) @>c>> \gme.
\end{CD}
\end{equation}
In local coordinates, this operator may be written as
$D=\sum\,c(dx^i)\,\n^\E_{\d_i}$. Note that $D$ sends even sections to odd
sections and vice versa: $D:\, \Gam(M,\E^\pm)\to \Gam(M,\E^\mp)$.

Suppose now that the Clifford module $\E$ is endowed with a Hermitian structure
and consider the $L^2$-scalar product on the space of sections $\gme$ defined by
the Riemannian metric on $M$ and the Hermitian structure on $\E$. By
\cite[Proposition~3.44]{BeGeVe}, {\em the Dirac operator associated to a
Clifford connection $\n^\E$ is formally self-adjoint with respect to this scalar
product if and only if $\E$ is a self-adjoint Clifford module and $\n^\E$ is a
Hermitian connection.}


Let $\L$ be a Hermitian line bundle endowed with a Hermitian connection $\n^\L$.
For each integer $k\ge0$, we consider the bundle $\Ek$ as a Clifford module with
Clifford action $c(a)\otimes1$. The tensor product connection $\n^{\Ek}$ on
$\Ek$ is a self-adjoint Clifford connection (cf. \refss{clmodule}) and, hence,
define a formally self-adjoint Dirac operator $D_k:\G(M,\E\otimes\Lk)\to
\G(M,\E\otimes\Lk)$.

We denote by $D_k^\pm$ the restriction of $D_k$ to the space
$\G(M,\E^\pm\otimes\Lk)$.

\ssec{Gamma}{A discrete group action}
Assume now that a discrete group $\G$ acts freely on $M$ by orientation
preserving isometries and that the quotient manifold $X=M/\G$ is compact. Then
$\G$ acts naturally on $C(M)$ preserving the fiberwise algebra structure.

Suppose that there are given actions of $\G$ on the bundles $\E,\L$ which cover
the action of $\G$ on $M$ and preserve the Hermitian structures and the
connections on $\E,\L$. We also assume that the Clifford action of $C(M)$ on $M$
is $\G$-invariant. Then the Dirac operator $D_k$ commutes with the $\G$-action.

It follows from \cite[\S3]{Atiyah76}, that $D_k$ extends to a self-adjoint
unbounded operator on the space $L^2(M,\Ek)$ of square-integrable sections of
$\Ek$.

\ssec{Dirvanish}{The vanishing theorem}
The curvature $F^\L$ of $\n^\L$ is an imaginary valued 2-form on $M$. If it is
non-degenerate at all points of $M$, then $iF^\L$ is a symplectic form on $M$,
and, hence, defines an orientation of $M$.

Our first result is the following
\th{Dirvanish}
Assume that the curvature $F^\L=(\n^\L)^2$ of the connection $\n^\L$ is
non-degenerate at all points of $M$. If the orientation defined by the
symplectic form $iF^\L$ coincides with the original orientation of $M$, then
\begin{equation}\label{E:Dirvanish}
        \Ker D^-_k \ = \ 0 \qquad \mbox{for} \qquad k\gg 0.
\end{equation}
Otherwise, $\Ker D^+_k = 0$ for $k\gg 0$.
\eth
The proof is given in \refss{prmain}. For the case when $\G$ is a trivial group
(and, hence, $M$ is a compact manifold) this theorem was established in
\cite{BrVanish}.
\rem{Dirvanish}
Since the operators $D_k^\pm$ are $\G$-invariant, they are lifts of certain
operator $D_{X,k}^\pm$ on $X$. By the $L^2$-index theorem of Atiyah,
\cite{Atiyah76}, the $\G$-index of $D_k^\pm$ is equal to the usual index of
$D_{X,k}^\pm$ (we refer the reader to \cite{Atiyah76} for the definitions of
$\G$-dimensions and $\G$-index). Combining with our vanishing theorem we obtain
\[
        \dim_\G\, \Ker D^\pm_k \ = \ \dim\, \Ker D_{X,k}^\pm
                \qquad \mbox{for any} \qquad k\gg 0.
\]
Here $\dim_\G$ denotes the $\G$-dimension, cf. \cite{Atiyah76}.
\erem
\rem{Griffiths}
\reft{Dirvanish} remains valid if $\L$ is a vector bundle of dimension higher
than 1 (cf. \cite{Griffiths65} for analogous generalization of the Kodaira
vanishing theorem). In this case $\Lk$ should be understood as the $k$-th
symmetric power of $\L$. Also the curvature $F^\L$ becomes a 2-form with values
in the bundle $\End(\L)$ of endomorphisms of $\L$. We say that it is
non-degenerate if, for any $0\not=\xi\in\L$, the imaginary valued form
$\<F^\L\xi,\xi\>$ is non-degenerate. In this case, the orientation of $M$
defined by the form $i\<F^\L\xi,\xi\>$ is independent of $\xi\not=0$. If this
orientation coincides with (resp. is opposite to) the original orientation of
$M$, then $\Ker D_k^-=0$ (resp. $\Ker{}D_k^+=0)$.

The proof is a combination of the methods of this paper with those of
\cite{Griffiths65}. The details will appear elsewhere.
\erem
\sec{cvanish}{Complex manifolds. The $L^2$ analogue of the Andreotti-Grauert theorem}

In this section we present an $L^2$ analogue of the Andreotti-Grauert vanishing
theorem. This result is a refinement of \reft{Dirvanish} for the case when $M$
is a complex manifold.

\ssec{L2coh}{The reduced $L^2$ Dolbeault cohomology}
Suppose $M$ is a complex manifold endowed with a holomorphic free action of a
discrete group $\G$, such that the quotient $X=M/\G$ is compact. Let $\W$ be a
holomorphic $\G$-equivariant vector bundle over $M$ and let $\L$ be a
holomorphic $\G$-equivariant line bundle over $M$.

Fix a $\G$-invariant Hermitian metric on $TM\otimes\CC$ and $\G$-invariant
Hermitian metrics on the bundles $\E,\L$. Let $\A^{0,*}(M,\E\otimes\L^k)$ and
$L^2\A^{0,*}(M,\E\otimes\L^k)$ denote respectively the spaces of smooth and
square-integrable $(0,*)$-differential forms on $M$ with values in
$\E\otimes\L^k$.

Set
\begin{align}
        Z^j \ &= \ \Ker\Big(\p:L^2\A^{0,j}(M,\E\otimes\L^k)\to
                  L^2\A^{0,j+1}(M,\E\otimes\L^k)\Big);\notag \\
        B^j \ &= \ \IM\Big(\p:L^2\A^{0,j-1}(M,\E\otimes\L^k)\to
                  L^2\A^{0,j}(M,\E\otimes\L^k)\Big)\notag
\end{align}
and let $\oB^j$ denote the closure of $B^j$ in $L^2\A^{0,j}(M,\E\otimes\L^k)$.

The {\em (reduced) $L^2$ Dolbeault cohomology} of $M$ with coefficients in the
bundle $\E\otimes\L^k$ is the quotient space
\[
        L^2H^j(M,\E\otimes\L^k) \ = \ Z^j/\oB^j.
\]
Note, that though the $L^2$ square product depends on the choices of Hermitian
metrics on $M,\E$ and $\L$, the topology of the Hilbert space
$L^2\A^{0,*}(M,\E\otimes\L^k)$ does not. So the cohomology
$L^2H^j(M,\E\otimes\L^k)$ is essentially independent of the metrics.

\rem{coh-dir}
The reduced $L^2$ cohomology is isomorphic to the kernel of the Dolbeault-Dirac
operator
$\sqrt2(\p+\p^*):L^2\A^{0,*}(M,\E\otimes\L^k)\to{}L^2\A^{0,*}(M,\E\otimes\L^k)$.
If $M$ is a \ka manifold, then (cf. \cite[Proposition~3.67]{BeGeVe}) the
Dolbeault-Dirac operator has the form \refe{dir1} (see \refss{acAG}). This
connects the material of this section with \reft{Dirvanish}. See \refr{2} for
more details.
\erem

\ssec{curv}{The curvature of the Chern connection}
Let $\n^\L$ be the {\em Chern connection \/} on $\L$, i.e., the unique
holomorphic connection which preserves the Hermitian metric. Then $\n^\L$ is
preserved by the action of $\G$. The curvature $F^\L$ of $\n^\L$ is a
$\G$-invariant $(1,1)$-form which is called the {\em curvature form of the
Hermitian metric $h^\L$}.

The orientation condition of \reft{Dirvanish} may be reformulated as follows.
Let $(z^1\nek z^n)$ be complex coordinates in the neighborhood of a point
$x\in{}M$. The curvature $F^\L$ may be written as
\[
        iF^\L \ = \ \frac{i}2\sum_{i,j} F_{ij}dz^i\wedge d\oz^j.
\]
Denote by $q$ the number of negative eigenvalues of the matrix $\{F_{ij}\}$.
Clearly, the number $q$ is independent of the choice of the coordinates. We will
refer to this number as the {\em number of negative eigenvalues of the curvature
$F^\L$ at the point $x$}. Then the orientation defined by the symplectic form
$iF^\L$ coincides with the complex orientation of $M$ if and only if $q$ is
even.

\ssec{AG}{The $L^2$ Andreotti-Grouert theorem}
A small variation of the method used in the proof of \reft{Dirvanish} allows to
get a more precise result which depends not only on the parity of $q$ but on $q$
itself. In this way we obtain the following
%
\th{AG}
Let $M$ be a complex manifold on which a discrete group $\G$ acts freely so that
$M/\G$ is compact. Let $\L$ be a $\G$-equivariant holomorphic line bundle over
$M$. Assume that $\L$ carries a $\G$-invariant Hermitian metric whose curvature
form $F^\L$ has at least $q$ negative and at least $p$ positive eigenvalues at
any point $x\in M$. Then, for any $\G$-equivariant holomorphic vector bundle
$\W$ over $M$, the cohomology $L^2H^{0,j}(M,\W\otimes\Lk))$ vanishes for
$j\not=q, q+1\nek n-p$ and $k\gg0$.
\eth
The proof is given in \refss{prAG}. If $\G$ is a trivial group, \reft{AG}
reduces to the classical Andreotti-Grauert vanishing theorem
\cite{AndGr62,DemPetSch93}.

Contrary to \reft{Dirvanish}, the curvature $F^\L$ in \reft{AG} needs not be
non-degenerate. If $F^\L$ is non-degenerate, then the number $q$ of negative
eigenvalues of $F^\L$ does not depend on the point $x\in M$. Then we obtain the
following
\cor{AG}
If, in the conditions of \reft{AG}, the curvature $F^\L$ is non-degenerate and
has exactly $q$ negative eigenvalues at any point $x\in M$, then
$L^2H^{0,j}(M,\W\otimes\Lk)$ vanishes for any $j\not=q$ and $k\gg0$.
\ecor
The most important is the case when the bundle $\L$ is positive, i.e., when the
matrix $\{F_{ij}\}$ is positive. In this case \refc{AG} generalizes the
classical Kodaira vanishing theorem (cf., for example,
\cite[Theorem~3.72(2)]{BeGeVe}) to covering manifolds.
\rem{2}
a. \ It is interesting to compare \refc{AG} with \reft{Dirvanish} for the case
when $M$ is a \ka manifold. In this case the Dirac operator $D_k$ is equal to
the Dolbeault-Dirac operator, cf. \cite[Proposition~3.67]{BeGeVe}. Hence (cf.
\refr{coh-dir}), \reft{Dirvanish} implies that $L^2H^{0,j}(M,\W\otimes\Lk)$
vanishes when the parity of $j$ is not equal to the parity of $q$. \refc{AG}
refines this result.

b. \ If $M$ is not a \ka manifold, then the Dirac operator $D_k$ defined by
\refe{dir1} is not equal to the Dolbeault-Dirac operator, and the kernel of
$D_k$ is not isomorphic to the cohomology $L^2H^{0,*}(M,\W\otimes\Lk)$. However,
we show in \refs{prAG} that the operators $D_k$ and $\sqrt2(\p+\p^*)$ have the
same asymptotic as $k\to\infty$. Then the vanishing of the kernel of $D_k$
implies the vanishing of the $L^2$-cohomology.
\erem
\rem{Griffiths2}
As in \refr{Griffiths}, one can generalize \reft{AG} to the case when $\L$ is a
vector bundle of dimension greater than 1. In this way one obtains, in
particular, an $L^2$ analogue of the Griffiths vanishing theorem
\cite{Griffiths65}.
\erem

\sec{almcomp}{Almost complex manifolds. An analogue of the Andreotti-Grauert theorem.}

In this section we refine \reft{Dirvanish} for the case when $M$ is an almost
complex manifold and $F^\L$ is a $(1,1)$-form. From another point of view, the
results of this section generalize the $L^2$ Andreotti-Grauert theorem to almost
complex manifolds, cf. \refr{UB}.b.

As an application we prove that {\em semiclassically the $\spin^c$-quantization
of an almost complex covering manifold gives an honest Hilbert space} (and not
just a virtual one).

\ssec{acAG0}{}
Let $M,\L,\G$ be as in \refss{Gamma}. Assume, in addition, that $M$ is an almost
complex $2n$-dimensional manifold, the action of $\G$ preserves the almost
complex structure $J$ on $M$, and the curvature $F^\L$ is a $(1,1)$-form on $M$
with respect to $J$. The later condition implies that, for any $x\in M$ and any
basis $(e^1\nek e^n)$ of the holomorphic cotangent space $(T\ha{M})^*$, one has
\[
        iF^\L \ = \ \frac{i}2\sum_{i,j} F_{ij}e^i\wedge \oe^j.
\]
We denote by $q$ the number of negative eigenvalues of the matrix $\{F_{ij}\}$.
As in \refss{curv}, the orientation of $M$ defined by the symplectic form
$iF^\L$ depends only on the parity of $q$. It coincides with the orientation
defined by $J$ if and only if $q$ is even.

We fix a Riemannian metric $g^{TM}$ on $M$, such that the almost complex
structure $J:TM\to TM$ is skew-adjoint with respect to $g^{TM}$.

\ssec{clstr}{A Clifford action on $\Lam^j(T^{0,1}M)^*$}
Let $\Lam^j =\Lam^j(T^{0,1}M)^*$ denote the bundle of $(0,j)$-forms on $M$ and
set
\[
        \Lam^+ \ = \ \bigoplus_{j \text{ even}} \Lam^j, \quad
        \Lam^-  \ = \ \bigoplus_{j \text{ odd}} \Lam^j.
\]

Let $\lam^{1/2}$ be the square root of the complex line bundle
$\lam=\det{}T\ha{M}$ and let $\calS$ be the spinor bundle over $M$ associated to
the Riemannian metric $g^{TM}$. Although $\lam^{1/2}$ and $\calS$ are defined
only locally, unless $M$ is a spin manifold, it is well known (cf.
\cite[Appendix~D]{LawMic89}) that the products $\calS^\pm\otimes\lam^{1/2}$ are
globally defined and $\Lam^\pm= \calS^\pm\otimes\lam^{1/2}$. Since the spinor
bundle is, by definition, a Clifford module, the last equality defines a
Clifford action of $C(M)$ on $\Lam$, cf. \refss{clmodule}. More explicitly, this
action may be described as follows: if $f\in\Gam(M,T^*M)$ decomposes as
$f=f\ha+f\ah$ with $f\ha\in\Gam(M,(T\ha{}M)^*)$ and $f\ah\in\Gam(M,(T\ah M)^*)$,
then the Clifford action of $f$ on $\alp\in\Gam(M,\Lam)$ equals
\begin{equation}\label{E:claction}
        c(f)\alp \ = \ \sqrt{2}\, \left( f\ah\wedge\alp  \ - \ \iot(f\ha)\, \alp\right).
\end{equation}
Here $\iot(f\ha)$ denotes the interior multiplication by the vector field
$(f\ha)^*\in T\ah M$ dual to the 1-form $f\ha$. This action is self-adjoint with
respect to the Hermitian structure on $\Lam$ defined by the Riemannian metric
$g^{TM}$ on $M$.

The Levi-Civita connection $\n^{TM}$ of $g^{TM}$ induces Hermitian connections
on $\lam^{1/2}$ and on $\calS$. Let
$\n^M={}\n^\calS\otimes1+1\otimes\n^{\lam^{1/2}}$ be the product connection (cf.
\refss{clmodule}). Then $\n^M$ is a well-defined Hermitian Clifford connection
on bundle $\Lam$.

Note also that the grading $\Lam=\Lam^+\oplus\Lam^-$ is natural.

\ssec{acAG}{An analogue of the Andreotti-Grauert theorem}
Let $\W$ be a $\G$-equivariant vector bundle over $M$ endowed with a
$\G$-invariant Hermitian metric and a $\G$-invariant Hermitian connection. Set
$\E^\pm=\Lam^\pm\otimes\W$. Then $\E=\E^+\oplus\E^-$ is a self-adjoint Clifford
module. Let $\n^\E$ be the product of the connection $\n^M$ on $\Lam$ and a
Hermitian connection on $\W$. Then $\n^\E$ is a Hermitian Clifford connection on
$\E$. The space $L^2(M,\E\otimes\Lk)$ of square-integrable sections of
$\E\otimes\Lk$ coincides with the space $L^2\A^{0,*}(M,\W\otimes\Lk)$ of
square-integrable differential forms of type $(0,*)$ with values in
$\W\otimes\Lk$. Let
\[
        D_k:\, L^2\A^{0,*}(M,\W\otimes\Lk) \ \to \ L^2\A^{0,*}(M,\W\otimes\Lk)
\]
denote the Dirac operator corresponding to the tensor product connection on
$\E\otimes\Lk$.

{}For a form $\alp\in L^2\A^{0,*}(M,\W\otimes\Lk)$, we denote by $\|\alp\|$ its
$L^2$-norm and by $\alp_i$ its component in $L^2\A^{0,i}(M,\W\otimes\Lk)$.
%
\th{UB+}
In the situation described above, assume that the matrix $\{F_{ij}\}$ has at
least $q$ negative and at least $p$ positive eigenvalues at any point $x\in M$.
Then there exists a sequence $\eps_1,\eps_2,\ldots$ convergent to zero, such that for
any $k\gg 0$ and any $\alp\in\Ker D_k^2$ one has
\[
        \|\alp_j\| \ \le \  \eps_k\|\alp\|, \qquad \mbox{for} \quad
                                                        j\not=q,q+1\nek n-p.
\]

In particular, if the form $F^\L$ is non-degenerate and $q$ is the number of
negative eigenvalues of $\{F_{ij}\}$ (which is independent of $x\in M$), then
there exists a sequence $\tileps_{1},\tileps_{2},\ldots$, convergent to zero,
such that $\alp\in\Ker D_k$ implies
\[
         \|\alp \ - \ \alp_q\| \ \le \ \tileps_{k}\|\alp_q\|.
\]
\eth
\reft{UB+} is proven in \refss{prUB}. For the case when $\G$ is a trivial group
it was proven in \cite{BrVanish}.

\rem{UB}
a. \ \reft{UB+} implies that, if $F^\L$ is non-degenerate, then $\Ker D_k$ is
dominated by the component of degree $q$. If $\alp\in L^2(M,\E^-)$ (resp.
$\alp\in L^2(M,\E^+)$) and $q$ is even (resp. odd) then $\alp_q=0$. So, we
obtain the vanishing result of \reft{Dirvanish} for the case when $M$ is almost
complex and $F^\L$ is a $(1,1)$-form.

b. \ \reft{UB+} is an analogue of \reft{AG}. Of course, the cohomology
$L^2H^{0,j}(M,\W\otimes\Lk)$ is not defined if $J$ is not integrable. Moreover,
the square $D^2_k$ of the Dirac operator does not preserve the $\ZZ$-grading on
$L^2\A^{0,*}(M,\W\otimes\Lk)$. Hence, one can not hope that the kernel of $D_k$
belongs to $\oplus_{j=q}^{n-p}L^2\A^{0,j}(M,\W\otimes\Lk)$. However, \reft{UB+}
shows, that for any $k\gg0$ and any $\alp\in\Ker D_k$, ``most of the norm" of
$\alp$ is concentrated in $\oplus_{j=q}^{n-p}L^2\A^{0,j}(M,\W\otimes\Lk)$.
\erem

\ssec{BU}{Positive line bundle. Quantization}
Probably the most interesting application of \reft{UB+} may be obtained by
choosing $\L$ to be a positive line bundle. Then \reft{UB+} (and, in fact, even
\reft{Dirvanish}) states that $\Ker D_k^- =0$. Assume also that $\W$ is a
trivial line bundle. In this case, the {\em index space}
$\Ker{}D_k^+\ominus\Ker{}D_k^-$ of $D_k$ plays the role of the
``quantum-mechanical space" (or the space of ``quantization") in the scheme of
geometric quantization (cf. \cite{Vergne92}). The geometric quantization for
$k\gg 0$ is called the {\em semiclassical limit}.

Unfortunately, the index space of $D_k$ is only a {\em virtual} (or graded)
vector space, in general. However, Theorems~\ref{T:Dirvanish} and \ref{T:UB+}
imply that {\em semiclassically this is an ``honest" vector space}. For the case
when $\G$ is the trivial group (so that $M$ is compact) this result was
established by Borthwick and Uribe \cite[Theorem~2.3]{BoUr96}.

Note, that \reft{UB+} implies also that, for large values of $k$, the quantum
mechanical-space $\Ker D_k^+$ is ``almost" a subspace of $\A^{0,0}(M,\Lk)$. More
precisely, the restriction of the projection
$L^2\A^{0,*}(M,\Lk)\to{}L^2\A^{0,0}(M,\Lk)$ to $\Ker D_k^+$ tends to the
identity operator when $k\to\infty$. This fact is very important in the study of
semiclassical properties of quantization, cf. \cite[\S4]{BoUr96}.

\ssec{estonD}{Estimate on the Dirac operator}
The main ingredient of the proof of \reft{UB+} (cf. \refss{prUB}) is the
following estimate on $D_k$, which also has an independent interest:
\prop{estDk}
If the matrix $\{F_{ij}\}$ has at least $q$ negative and at least $p$ positive
eigenvalues at any point $x\in M$, then there exists a constant $C>0$, such that
\[
        \|D_k\, \alp\| \ \ge \ Ck^{1/2}\, \|\alp\|,
\]
for any $k\gg0, \ j\not=q,q+1\nek n-p$ and $\alp\in
L^2\A^{0,j}(M,\W\otimes\Lk)\cap \A^{0,j}(M,\W\otimes\Lk)$.
\eprop
The proof is given in \refss{prestDk}.

\sec{repr}{Application to the representation theory}

In this section we explain very briefly how Theorems~\ref{T:Dirvanish} and
\ref{T:AG} can be applied to the study of homogeneous vector bundles. In
particular we recover a vanishing theorem which was originally conjectured by
Langlands \cite{Langlands66} and proven by Griffiths and Schmid
\cite{GrifSch69}. We also indicate how one can use certain generalizations of
\reft{Dirvanish} to get a new proof of the vanishing theorem of Atiyah and
Schmid \cite[Th.~5.20]{AtSch77}. The details will appear in a separate paper.

\ssec{repr}{New proof of a theorem of Griffiths and Schmid}
Let $G$ be a connected non-compact real semi-simple Lie group and assume that it
has a compact Cartan subgroup $H$. Let $K\supset H$ be a maximal compact
subgroup of $G$. Let $\grg,\grk,\grh$ denote the Lie algebras of $G,K,H$ and let
$\grg_\CC,\grk_\CC,\grh_\CC$ denote their complexifications.

Denote by $\Del$ the set of roots for $(\grg_\CC, \grh_\CC)$. It decomposes as a
disjoint union $\Del=\Del_c\cup\Del_n$, where $\Del_c$ and $\Del_n$ are
respectively the set of compact and non-compact roots with respect to $K\subset
G$. Choose a system $P\subset\Del$ of positive roots.

Let $M=G/H$ and let $\L_\lam\to M$ be the line bundle over $M$ induced by the
character $\lam$ of $H$. It is well known (cf. \cite[\S1]{GrifSch69}) that the
choice of positive root system $P$ defines a complex structure on $M$ and that
$\L_\lam$ has a natural structure of a holomorphic line bundle over $M$.

By a theorem of Borel \cite{Borel62}, there exists a discrete subgroup
$\G\subset G$ which acts freely on $M$ and such that the quotient space
$\G\backslash{}M=\G\backslash{}G/H$ is compact. This allows us to apply
\reft{AG} to the study of $L^2$ cohomology of $\L_\lam$. As a result we obtain a
new proof of the following theorem, which was originally conjectured by
Langlands \cite{Langlands66} and proven by Griffiths and Schmid
\cite[Th.~7.8]{GrifSch69}.

\th{GrSch}
Let $(\cdot,\cdot)$ denote the scalar product on $\grh^*$ induced by the
Cartan-Killing form on $\grh$. Suppose $\lam$ is a character of $H$, such that
$(\lam,\alp)\not=0$ for any $\alp\in\Del$ and set
\[
        \iot(\lam) \ = \ \#\{ \alp\in P\cap\Del_c: \ (\lam, \alp) <0\, \}
                \ + \ \#\{ \alp\in P\cap\Del_n: \ (\lam, \alp) >0\, \}.
\]
Let $\rho$ denotes the half-sum of the positive roots of $(\grg,\grh)$. There
exists an integer $m$ such that, if $k\ge m$, then
\[
        L^2H^{0,j}(M,\calL_{k\lam}) \ = \ 0
                \quad \mbox{for} \quad j\not=\iot(\lam+\rho).
\]
\eth
\prf
By \cite[Th.~4.17D]{GrifSch69}, the bundle $\calL$ possesses a Chern connection
whose curvature is non-degenerate and has exactly $\iot(\lam+\rho)$ negative
eigenvalue. \reft{GrSch} follows now from \refc{AG}.
\eprf
\rem{GrSch}
1. \ The significance of \reft{GrSch} is that it allows to prove (cf.
\cite{GrifSch69}) that, for $j=\iot(\lam+\rho), k\gg0$, the space
$L^2H^{0,j}(M,\calL_{k\lam})$ is an irreducible discrete series representation
of $G$. In this way ``most" of the discrete series representation may be
obtained.

2. \ \reft{GrSch} may be considerably improved. In particular (cf.
\cite[Th.~7.8]{GrifSch69}) there exists a constant $b$, depending only on $G$
and $H$, such that
\[
        L^2H^{0,j}(M,\calL_{\lam}) \ = \ 0
                \quad \mbox{for} \quad j\not=\iot(\lam+\rho).
\]
whenever $\lam$ satisfies $|(\lam,\alp)|>b$ for every $\alp\in\Del$.

3. \ Let $G_\CC$ be a complex form of $G$ and suppose that $B\subset H$ is a
Borel subgroup of $G_\CC$, such that $V=G\cap B$ is compact. Let $L_\lam$ denote
the irreducible $V$ module with highest weight $\lam$ and let
$\calL_\lam=G\times_V L_\lam$ be the corresponding homogeneous vector bundle
over $M=G/V$. Using the generalization of \reft{AG} discussed in
\refr{Griffiths2}, one can show that \reft{GrSch} remains true in this case. In
this form the theorem is proven in \cite{GrifSch69}.
\erem

\ssec{AtSch}{A vanishing theorem of Atiyah and Schmid}
Let $L_\lam$ denote the irreducible $K$ module with highest weight $\lam$ and
let $\calL_\lam=G\times_K L_\lam$ be the corresponding homogeneous vector bundle
over $M=G/K$. For simplicity, assume also that $M$ possesses a $G$-equivariant
spinor bundle $\calS=\calS^+\oplus\calS^-$. Then (cf. \cite{AtSch77}) there is a
canonically defined Dirac operator
$D_\lam^\pm:L^2(M,\L_\lam\otimes\calS^\pm)\to{}L^2(M,\L_\lam\otimes\calS^\mp)$.
Using the generalization of \reft{Dirvanish} discussed in \refr{Griffiths}, one
can show that $\Ker D_\lam=0$ if $\lam$ is {\em sufficiently non-singular} (i.e.
if $|(\lam,\alp)|\gg0$ for any $\alp\in\Del$). This is a particular case of
\cite[Th.~5.20]{AtSch77}.


\sec{prdirvanish}{Proof of the vanishing theorem for the half-kernel of a Dirac operator}

In this section we present a proof of \reft{Dirvanish} based on
Propositions~\ref{P:D-D} and \ref{P:lapl}, which will be proved in the following
sections.

The idea of the proof is to study the large $k$ behavior of the square $D_k^2$
of the Dirac operator.

\ssec{tJ}{The operator $\tJ$}
We need some additional definitions. Recall that $F^\L$ denotes the curvature of
the connection $\n^\L$. In this subsection we do not assume that $F^\L$ is
non-degenerate. For $x\in M$, define the skew-symmetric linear map
$\tJ_x:T_xM\to T_xM$ by the formula
\[
        iF^\L(v,w) \ = \ g^{TM}(v,\tJ_xw), \qquad v,w\in T_xM.
\]
The eigenvalues of $\tJ_x$ are purely imaginary.
Define
\begin{equation}\label{E:tau}
        \tau(x) \ = \ \Tr^+ \tJ_x \ := \ \mu_1+\cdots+\mu_l, \qquad
        m(x)    \ = \ \min_{j}\, \mu_j(x).
\end{equation}
where $i\mu_j, \ j=1\nek l$ are the eigenvalues of $\tJ_x$ for which $\mu_j>0$.
Note that $m(x)$ is well defined and positive if the curvature $F^\L$ does not
vanishes at the point $x\in M$.

\ssec{lapl}{Estimate on $D_k^2$}
Consider the {\em rough (or metric) Laplacian}
\[
        \Del_k \ := \ (\n^{\E\otimes\Lk})^*\, \n^{\E\otimes\Lk},
\]
where $(\n^{\Ek})^*$ denote the formal adjoint of the covariant derivative
\[
        \n^\Ek:C^\infty_0(M,\Ek){\to}C^\infty_0(M,\Ek\otimes T^*M).
\]
Here, as usual, $C^\infty_0$ denotes the space of sections with compact support.

Since $\Del_k$ is an elliptic $\G$-invariant operator, it follows from
\cite[Proposition~3.1]{Atiyah76} that it is a self-adjoint operator.

Our estimate on the square $D_k^2$ of the Dirac operator is obtained in two
steps: first we compare it to the rough Laplacian $\Del_k$ and then we estimate
the large $k$ behavior of $\Del_k$. These two steps are the subject of the
following two propositions.
%
\prop{D-D}
a. \ For any integer $k$, the difference $D_k^2-\Del_k$ is a bounded operator on
$L^2(M,\Ek)$.

b. \ Supposed that the differential form $F^\L$ is non-degenerate. If the
orientation defined on $M$ by the symplectic form $iF^\L$ coincides with (resp.
is opposite to) the given orientation of $M$, then there exists a constant $C$
such that, for any $s\in L^2(M,\E^-\otimes\Lk)$ (resp. for any
$s{\in}L^2(M,\E^+\otimes\Lk)$), one has an estimate
\[
        \big\langle\, (D^2_k - \Del_k)\, s,\, s\, \big\rangle \ \ge \
                -k\,\big\langle\, (\tau(x)-2m(x))\, s,\, s\, \big\rangle - C\, \|s\|^2.
\]
Here $\<\cdot,\cdot\>$ denotes the $L^2$-scalar product on the space of sections
and $\|\cdot\|$ is the norm corresponding to this scalar product.
\eprop
The proposition is proven in \refss{lichn} using the Lichnerowicz formula \refe{lichn}.

Set
\[
        \Dom(D_k^2) \ := \
                \left\{ s\in L^2(M,\Ek):  \ D_k^2s\in L^2(M,\Ek) \right\}.
\]
where $D_k^2s$ is understood in the sense of distributions. In the next
proposition we do not assume that $F^\L$ is non-degenerate.
\prop{lapl}
Suppose that $F^\L$ does not vanish at any point $x\in M$. For any $\eps>0$,
there exists a constant $C_\eps$ such that, for any $k\in\ZZ$ and any
$s{\in}\Dom(D_k^2)$ one has
\eq{lapl}
        \<\Del_k\, s,s\> \ \ge \ k\, \<(\tau(x)-\eps) s,s\> \ - \ C_\eps\, \|s\|^2.
\end{equation}
\eprop
\refp{lapl} is proven in \refs{lapl}. For the case when $\G$ is a trivial group
(so that $M$ is compact) it was essentially proven in \cite[Theorem~2.1]{BoUr96}
(see also \cite[Proposition~4.4]{BrVanish}).

\ssec{prmain}{Proof of \reft{Dirvanish}}
Assume that the orientation defined by $iF^\L$ coincides with the given
orientation of $M$ and $s\in L^2(M,\E^-\otimes\L)$, or that the orientation
defined by $iF^\L$ is opposite to the given orientation of $M$ and
$s{\in}L^2(M,\E^+\otimes\L)$. By \refp{D-D},
\eq{D>D}
        \<D_k^2\, s,s\> \ \ge \ \< \Del_k\, s,s\> \ - \
        k\,\big\langle\, (\tau(x)-2m(x))\, s,\, s\big\rangle \ - \ C\, \|s\|^2.
\end{equation}
Choose
\[
        0 \ < \ \eps \ < \ 2\, \min_{x\in M} m(x)
\]
and set
\[
        C'=2\min_{x\in M} m(x)-\eps>0.
\]
Assume now that $D_k^2s\in L^2(M,\Ek)$. Then, it follows from \refp{D-D}.a, that
$\Del_ks\in L^2(M,\Ek)$. Using \refe{lapl} and \refe{D>D}, we obtain
\begin{equation}\label{E:Dk2>}
         \<D_k^2\, s,\, s\> \ \ge \ \big(k C'-(C+C_\eps)\big)\, \|s\|^2.
\end{equation}
Thus, for $k> (C+C_\eps)/C'$, we have $\<D_k^2\, s,s\>>0$. Hence, $D_k s\not=0$.
\hfill$\square$

\sec{prUB}{Proof of the Andreotti-Grauert-type theorem for almost complex manifolds}

In this section we prove \reft{UB+} and \refp{estDk}. The proof is very similar
to the proof of \reft{Dirvanish} (cf. \refs{prdirvanish}). It is based on
\refp{lapl} and the following refinement of \refp{D-D}:

\prop{D-Dac}
Assume that the matrix $\{F_{ij}\}$ (cf. \refss{curv}) has at least $q$ negative
eigenvalues at any point $x\in M$. For any $x\in M$, we denote by $m_q(x)>0$ the
minimal positive number, such that at least $q$ of the eigenvalues of
$\{F_{ij}\}$ do not exceed $-m_q$. Then there exists a constant
 $C$  such that
\[
        \big\langle\, (D^2_k - \Del_k)\, \alp,\, \, \alp\big\rangle \ \ge \
                -k\, \big\langle\, (\tau(x)-2m_q(x))\, \alp,\, \alp\, \big\rangle
                \  - \ C\, \|\alp\|^2
\]
for any $j=0\nek q-1$ and any $\alp\in L^2\A^{0,j}(M,\W\otimes\Lk)$.
\eprop
The proposition is proven in \refss{prD-Dac}.

\ssec{prestDk}{Proof of \refp{estDk}}
Choose $0<\eps<2\min_{x\in M}\, m_q(x)$ and set
\[
        C' \ = \ 2\min_{x\in M}\, m_q(x)-\eps.
\]
{}Fix $j=0\nek q-1$ and $\alp\in L^2\A^{0,j}(M,\W\otimes\Lk)$ such that
$D_k^2\alp\in L^2\A^{0,*}(M,\W\otimes\Lk)$. By \refp{D-D}.a,
$\Del_k\alp{\in}L^2\A^{0,*}(M,\W\otimes\Lk)$. It follows from
Propositions~\ref{P:lapl} and \ref{P:D-Dac}, that
\[
        \< D_k^2\, \alp,\alp\> \ \ge \ kC'\, \|\alp\|^2 \ - \ (C+C_\eps)\, \|\alp\|^2.
\]
Hence, for any $k>2(C+C_\eps)/C'$, we have
\[
        \|D_k\, \alp\|^2 \ = \ \<D_k^2\, \alp,\, \alp\, \> \ \ge \
                \frac{k C'}2\, \|\alp\|^2.
\]
This proves \refp{estDk} for $j=0\nek q-1$. The statement for $j=n-p+1\nek n$
may be proven by a verbatim repetition of the above arguments, using a natural
analogue of \refp{D-Dac}. (Alternatively, the statement for $j=n-p+1\nek n$ may
be obtained as a formal consequence of the statement for $j=0\nek q-1$ by
considering $M$ with an opposite almost complex structure).
\hfill$\square$

\ssec{gr}{}
If the manifold $M$ is not K\"ahler, then the operator $D_k^2$ does not preserve
the $\ZZ$-grading on $\A^{0,*}(M,\W\otimes\Lk)$. However, the next proposition
shows that the {\em mixed degree operator} $\alp_i\mapsto(D_k^2\alp_i)_j$ is, in
a certain sense, small.
\prop{gr}
Set $\Dom^i(D_k^2)=
\{\alp\in L^2\A^{0,i}(M,\W\otimes\Lk):\, D_k^2\alp\in L^2\A^{0,*}(M,\W\otimes\Lk)\}$,
where $D_k^2\alp$ is understood in the sense of distributions. There exists a
sequence $\del_{1}, \del_{2}\nek$ such that $\lim_{k\to\infty}\del_{k}=0$ and
\[
        |\<D_k^2\, \alp, \bet\>| \ \le \del_{k}\, \<D_k^2\, \alp, \alp\> \ + \
          \del_{k}\, \<D_k^2\, \bet, \bet\> \ + \
             \del_{k}k\, \|\alp\|^{2} \ + \  \del_{k}k\, \|\bet\|^{2},
\]
for any $i\not=j$ and any $\alp\in \Dom^i(D_k^2), \bet\in \Dom^j(D_k^2)$.
\eprop
The proof of the proposition, based on the Lichnerowicz formula and
\cite{Dem85}, is given in \refss{prgr}.

\ssec{prUB}{Proof of \reft{UB+}}
Let $\alp\in \Ker D_k$ and fix $j\not\in q,q+1\nek n-p$. Set $\bet=\alp-\alp_j$.
Then
\[
        0 \ = \ \|D_k\, \alp\|^2 \ = \ \|D_k\, \alp_j\|^2  \ + \
               2\RE\, \<\, D_k\, \alp_j,\, D_k\, \bet\, \> \ + \
                  \|D_k\, \bet\|^2.
\]
Hence, it follows from \refp{gr}, that
\begin{equation}\label{E:alpbet}
        0 \ \ge \ (1-2\del_k)\, \|D_k\, \alp_j\|^2 \ + \
           (1-2\del_k)\, \|D_k\, \bet\|^2 \ - \
                  2\del_{k}k\, \|\alp_i\|^{2} \ - \ 2\del_{k}k\, \|\bet\|^{2}.
\end{equation}
If we assume now that $k$ is large enough, so that $1-2\del_k>0$, then we obtain
from \refe{alpbet} and \refp{estDk} that
\[
        \big((1-2\del_k)C^2k -2\del_kk\big)\,  \|\alp_i\|^{2} \ \le \
                        2\del_{k}k\, \|\bet\|^{2} \ \le \ 2\del_{k}k\, \|\alp\|^{2}.
\]
Thus
\[
        \|\alp_i\|^{2} \ \le \
                \frac{2\del_{k}}{(1-2\del_k)C^2 -2\del_k}\, \|\alp\|^{2}.
\]
Hence, \reft{UB+} holds with
$\eps_k=\sqrt{\frac{2\del_{k}}{(1-2\del_k)C^2-2\del_k}}$.
\hfill$\square$

\sec{prAG}{Proof of the $L^2$ Andreotti-Grauert theorem}

In this section we use the results of \refs{almcomp} in order to prove our $L^2$
version of the Andreotti-Grauert (\reft{AG}).

Note, first, that, if the manifold $M$ is K\"ahler, then the $L^2$
Andreotti-Grauert theorem follows directly from \reft{UB+}. Indeed, in this
case, the Dirac operator $D_k$ is equal to the Dolbeault-Dirac operator
$\sqrt2(\p+\p^*)$, cf. \cite[Proposition~3.67]{BeGeVe}. Hence, the restriction
of the kernel of $D_k$ to $L^2\A^{0,j}(M,\O(\W\otimes\Lk))$ is isomorphic to the
reduced $L^2$-cohomology $L^2H^j(M,\O(\W\otimes\Lk))$.

In general, $D_k\not=\sqrt2(\p+\p^*)$. However, the following proposition shows
that those two operators have the same ``large $k$ behavior".

\prop{D-dd}
Set $\E=\Lam(T^{0,1}M)\otimes\W$. In the conditions of \reft{AG}, there exists a
$\G$-invariant bundle map $A\in\End(\E)\subset{}\End(\E\otimes\Lk)$, independent
of $k$, such that
\begin{equation}\label{E:D-dd}
        \sqrt2\, \left(\p+\p^*\right) \ = \ D_k \ + \ A.
\end{equation}
\eprop
\prf
Choose a holomorphic section $e(x)$ of $\L$ over an open set $U\subset M$. It
defines a section $e^k(x)$ of $\Lk$ over $U$ and, hence, a holomorphic
trivialization
\begin{equation}\label{E:triv}
        U\times\CC \ \overset{\sim}{\longrightarrow} \
                \Lk, \qquad (x,\phi)\mapsto \phi\cdot e^k(x)\in \Lk
\end{equation}
of the bundle $\Lk$ over $U$. Similarly, the bundles $\W$ and $\W\otimes\Lk$ may
be identified over $U$ by the formula
\begin{equation}\label{E:WLk=W}
        w \ \mapsto \ w\otimes e^k.
\end{equation}

Let $h^\L$ and $h^\W$ denote the Hermitian fiberwise metrics on the bundles $\L$
 and  $\W$  respectively. Let  $h^{\W\otimes\Lk}$  denote the
Hermitian metric on $\W\otimes\Lk$ induced by the metrics $h^\L, h^\W$. Set
\[
        f(x)\ := \ |e(x)|^2, \qquad x\in U,
\]
where $|\cdot|$ denotes the norm defined by the metric $h^\L$. Under the
isomorphism \refe{WLk=W} the metric $h^{\W\otimes\Lk}$ corresponds to the metric
\begin{equation}\label{E:hk}
        h_k(\cdot,\cdot) \ = \ f^k\, h^\W(\cdot,\cdot)
\end{equation}
on $\W$.

By \cite[p.~137]{BeGeVe}, the connection $\n^\L$ on $\L$ corresponds under the
trivialization \refe{triv} to the operator
\[
        \Gam(U,\CC) \ \to \ \Gam(U,T^*U\otimes\CC);     \qquad
                        s\mapsto ds+kf^{-1}\d f\wedge s.
\]
Similarly, the connection on $\E\otimes\Lk= \Lam(T\ah M)^*\otimes\W\otimes\Lk$
corresponds under the isomorphism \refe{WLk=W} to the connection
\[
        \n_k:\, \alp \ \mapsto \ \n^\E\alp \ + \ kf^{-1}\d f\wedge\alp, \qquad
                \alp\in \Gam(U,\Lam(T\ah U)^*\otimes\W|_U)
\]
on $\E|_U$. It follows now from \refe{claction} and \refe{dir1} that the Dirac
operator $D_k$ corresponds under \refe{WLk=W} to the operator
\begin{equation}\label{E:tilDk}
        \tilD_k:\, \alp \ \mapsto \ D_0\alp \ - \ \sqrt2kf^{-1}\iot(\d f)\alp, \qquad
                        \alp\in \A^{0,*}(U,\W|_U).
\end{equation}
Here $\iot(\d f)$ denotes the contraction with the vector field $(\d f)^*\in
T\ah{M}$ dual to the 1-form $\d{f}$, and $D_0$ stands for the Dirac operator on
the bundle $\E=\E\otimes\L^0$.

Let $\p_k^*:\A^{0,*}(U,\W|_U) \to \A^{0,*-1}(U,\W|_U)$ denote the formal adjoint
of the operator $\p$ with respect to the scalar product on $\A^{0,*}(U,\W|_U)$
determined by the Hermitian metric $h_k$ on $\W$ and the Riemannian metric on
$M$. Then, it follows from \refe{hk}, that
\begin{equation}\label{E:pk}
        \p_k^* \ = \ \p \ + \ kf^{-1}\iot(\d f).
\end{equation}
By \refe{tilDk} and \refe{pk}, we obtain
\[
        \sqrt2\, \left(\p+\p^*_k\right) \ - \ \tilD_k \ =  \
                \sqrt2\, \left(\p+\p^*_0\right) \ - \ D_0.
\]
Set $A=\sqrt2(\p+\p^*_0)-D_0$. By \cite[Lemma~5.5]{Duis96}, $A$ is a zero order
operator, i.e. $A\in\End(\E)$ (note that our definition of the Clifford action
on $\Lam(T\ah{M})^*$ and, hence, of the Dirac operator defers from \cite{Duis96}
by a factor of $\sqrt2$).

Since both operators $D_0$ and $\sqrt2(\p+\p^*_0)-D_0$ are $\G$-invariant, so is
$A$.
\eprf

\ssec{prAG}{Proof of \reft{AG}}
Since the bundle map $A\in\End(\E)$ defined in \refp{D-dd} is $\G$-invariant, it
defines a bounded operator
$L^2\A^{0,*}(M,\E\otimes\Lk)\to{}L^2\A^{0,*}(M,\E\otimes\Lk)$. Let
\[
        \|A\| \ = \ \sup_{\|\alp\|=1} \, \|A\alp\|, \qquad
                        \alp\in L^2\A^{0,*}(M,\E\otimes\Lk)
\]
be the norm of this operator. By \refp{estDk}, there exists a constant $C>0$
such that
\[
        \|D_k\, \alp\| \ \ge \ Ck^{1/2}\, \|\alp\|,
\]
for any $k\gg0, \ j\not=q,q+1\nek n-p$ and
$\alp\in{}L^2\A^{0,j}(M,\W\otimes\Lk)\cap\A^{0,j}(M,\W\otimes\Lk)$. Then, if
$k>\|A\|^2/C^2$, we have
\[
        \|\sqrt2(\p+\p^*)\alp\| \ = \ \|(D_k+A)\alp\| \ \ge \
        \|D_k\alp\| \ - \ \|A\|\, \|\alp\|
        \ge \ \left(Ck^{1/2}-\|A\|\right)\, \|\alp\| \ > \ 0,
\]
for any $j\not=q,q+1\nek n-p$ and $0\not=\alp\in{}L^2\A^{0,j}(M,\W\otimes\Lk)$.
Hence, the restriction of the kernel of the Dolbeault-Dirac operator to the
space $L^2\A^{0,j}(M,\E\otimes\Lk)$ vanishes for $j\not=q,q+1\nek n-p$.
\hfill$\square$

\sec{compare}{The Lichnerowicz formula. Proof of Propositions~\ref{P:D-D},
         \ref{P:D-Dac} and \ref{P:gr}}

In this section we use the Lichnerowicz formula (cf. \refss{lichn}) to prove the
Propositions~\ref{P:D-D}, \ref{P:D-Dac} and \ref{P:gr}.

\ssec{curv2}{The curvature of a Clifford connection}
Before formulating the Lichnerowicz formula, we need some more information about
Clifford modules and Clifford connections (cf. \cite[Section~3.3]{BeGeVe}).

Let $\n^\E$ be a Clifford connection on a Clifford module $\E$ and let
$F^\E=(\n^\E)^2\in\A^2(M,\End(\E))$ denote the curvature of $\n^\E$.

Let $\EndC(\E)$ denote the bundle of endomorphisms of $\E$ commuting with the
action of the Clifford bundle $C(M)$. Then the bundle $\End(\E)$ of all
endomorphisms of $\E$ is naturally isomorphic to the tensor product
\begin{equation}\label{E:End}
        \End(\E) \ \cong \ C(M)\otimes \EndC(\E).
\end{equation}

By Proposition~3.43 of \cite{BeGeVe}, $F^\E$ decomposes with respect to
\refe{End} as
\begin{equation}\label{E:curv}
        F^\E \ = \ R^\E \ + \ \fes, \qquad\qquad R^\E\in \A^2(M,C(M)), \
                                                \fes\in \A^2(M,\EndC(\E)).
\end{equation}
In this formula, $\fes$ is an invariant of $\n^\E$ called the {\em twisting
curvature \/} of $\E$, and $R^\E$ is determined by the Riemannian curvature $R$
of
 $M$. If  $(e_1\nek e_{2n})$  is an orthonormal frame of the tangent space
$T_xM, \ x\in M$ and $(e^1\nek e^{2n})$ is the dual frame of the cotangent space
 $T^*M$, then
\[
        R^\E(e_i,e_j) \ = \
          \frac14\, \sum_{k,l}\,
                \big\langle R(e_i,e_j)e_k,e_l\big\rangle \, c(e^k)\, c(e^l).
\]

\rem{twist}
Assume that $\calS$ is a spinor bundle (\cite[\S3.3]{BeGeVe}),
$\E=\calW\otimes\calS$ and $\n^\E$ is given by the tensor product of a Hermitian
connection on $\W$ and the Levi-Civita connection on $\calS$. Then
$\A(M,\EndC(\E))\cong\A(M,\End(\calW))$ and the twisting curvature $\fes$ is
equal to the curvature $F^\calW=(\n^\calW)^2$ via this isomorphism (cf.
\cite[p.~121]{BeGeVe}). This explains why $\fes$ is called the twisting
curvature.
\erem

\ssec{lichn}{The Lichnerowicz formula}
Let $\E$ be a Clifford module endowed with a Hermitian structure and let
$D:\gme\to \gme$ be a self-adjoint Dirac operator associated to a Hermitian
Clifford connection $\n^\E$. Consider the rough Laplacian (cf. \refss{lapl})
\[
        \Del^\E \ = \ (\n^\E)^*\, \n^\E\, : \ \gme \ \to \gme,
\]
where $(\n^\E)^*$ denotes the formal adjoint of
$\n^\E:\gme\to\Gam(M,T^*M\otimes\E)$. By \cite[Proposition~3.1]{Atiyah76}, the
operator $\Del^\E$ is self-adjoint.

The following {\em Lichnerowicz formula \/} (cf. \cite[Theorem~3.52]{BeGeVe})
plays a crucial role in our proof of vanishing theorems:
\begin{equation}\label{E:lichn}
        D^2 \ = \ \Del^\E \ + \ \bfc(\fes) \ + \ \frac{r_M}4.
\end{equation}
Here $r_M$ stands for the scalar curvature of $M$, $\fes$ is the twisting
curvature of $\n^\E$ and
\begin{equation}\label{E:bfc}
        \bfc(F) \ := \ \sum_{i<j}\, F(e_i,e_j)\, c(e^i)\, c(e^j), \qquad
                        F\in \A^2(M,\End(\E)),
\end{equation}
where $(e_1\nek e_{2n})$ is an orthonormal frame of the tangent space to $M$,
and $(e^1\nek e^{2n})$ is the dual frame of the cotangent space.

Let $\L$ be a Hermitian line bundle over $M$ endowed with a Hermitian connection
$\n^\L$ and let $\n_k= \n^{\E\otimes\Lk}$ denote the product connection (cf.
\refe{nEW}) on $\E\otimes\Lk$. It is a Hermitian Clifford connection on
$\E\otimes\Lk$. The twisting curvature of $\n_k$ is given by
\begin{equation}\label{E:Ek/S}
        F^{(\Ek)/\calS} \ = \ kF^\L \ + \ F^{\E/\calS}.
\end{equation}

We denote by $D_k$ and $\Del_k$ the Dirac operator and the rough Laplacian
associated to this connection. By \refe{Ek/S}, it follows from the Lichnerowicz
formula \refe{lichn}, that
\begin{equation}\label{E:lechnW}
        D_k^2 \ = \ \Del_k \ + \ k\, \bfc(F^\L) \ + \ A,
\end{equation}
where $F^\L=(\n^\L)^2$ is the curvature of $\n^\L$ and
\begin{equation}\label{E:A}
        A \ := \ \bfc(\fes) \ + \ \frac{r_M}4 \ \in \ \End(\E) \ \subset \
                \End(\E\otimes\Lk)
\end{equation}
is independent of $\L$ and $k$.

\ssec{cFL}{Calculation of  $\bfc(F^\L)$}
To compare $D_k^2$ with the Laplacian $\Del_k$ we now need to calculate the
operator $\bfc(F^\L)\in \End(\E)\subset \End(\E\otimes\Lk)$. This may be
reformulated as the following problem of linear algebra.

Let $V$ be an oriented Euclidean vector space of real dimension $2n$ and let
$V^*$ denote the dual vector space. We denote by $C(V)$ the Clifford algebra of
$V^*$. Let $E$ be a module over $C(V)$. We will assume that $E$ is endowed with
a Hermitian scalar product such that the operator $c(v):E\to E$ is
skew-symmetric for any $v\in V^*$. In this case we say that $E$ is a {\em
self-adjoint \/} Clifford module over $V$.

The space $E$ possesses a {\em natural grading \/} $E=E^+\oplus E^-$, where
$E^+$ and $E^-$ are the eigenspaces of the chirality operator with eigenvalues
$+1$ and $-1$ respectively, cf. \refss{chir}.

In our applications $V$ is the tangent space $T_xM$ to $M$ at a point $x\in M$
 and  $E$  is the fiber of  $\E$  over  $x$.

Let $F$ be an imaginary valued antisymmetric bilinear form on $V$. Then $F$ may
be considered as an element of $V^*\wedge V^*$. We need to estimate the operator
$\bfc(F)\in\End(E)$. Here $\bfc:\Lam V^*\to C(V)$ defined exactly as in
\refe{bfc} (cf. \cite[\S{3.1}]{BeGeVe}).

Let us define the skew-symmetric linear map $\tJ:V\to V$ by the formula
\[
        i F(v,w) \ = \ \<v,\tJ w\>, \qquad v,w\in V.
\]
The eigenvalues of $\tJ$ are purely imaginary. Let $\mu_1\ge\cdots\ge\mu_l>0$ be
the positive numbers such that $\pm{i}\mu_1\nek\pm{i}\mu_l$ are all the non-zero
eigenvalues of $\tJ$. Set
\[
        \tau \ = \ \Tr^+ \tJ\ := \ \mu_1+\cdots+\mu_l, \qquad
        m    \ = \ \min_{j}\, \mu_j.
\]

Clearly, in the conditions of \refp{D-D}, the bundle map $A$, defined in
\refe{A}, is invariant with respect to $\G$. Thus it defines a bounded operator
$L^2(M,\Ek)\to L^2(M,\Ek)$. Hence, by the Lichnerowicz formula \refe{lichn},
\refp{D-D} is equivalent to the following
\prop{cFL}
Suppose that the bilinear form $F$ is non-degenerate. Then it defines an
orientation of $V$. If this orientation coincides with (resp. is opposite to)
the given orientation of $V$, then the restriction of $\bfc(F)$ onto $E^-$
(resp. $E^+$) is greater than $-(\tau-2m)$, i.e., for any $\alp\in E^-$ (resp.
$\alp\in E^+)$
\[
        \<c(F)\alp,\alp\> \ \ge \ -(\tau-2m)\, \|\alp\|^2.
\]
\eprop
We will prove the proposition in \refss{prcFL} after introducing some additional
constructions. Since we need these constructions also for the proof of
\refp{D-Dac}, we do not assume that $F$ is non-degenerate unless this is stated
explicitly.

\ssec{cs}{A choice of a complex structure on  $V$}
By the Darboux theorem (cf. \cite[Theorem~1.3.2]{Audin91}), one can choose an
orthonormal basis $f^1\nek f^{2n}$ of $V^*$, which defines the positive
orientation of $V$ (i.e., $f^1\wedge\dots\wedge f^{2n}$ is a positive volume form
on $V$) and such that
\begin{equation}\label{E:iF}
       iF_x^\L \ = \  \sum_{j=1}^l\, r_j\, f^j\wedge f^{j+n},
\end{equation}
for some integer $l\le n$ and some non-zero real numbers $r_j$. We can and we
will assume that $|r_1|\ge|r_2|\ge\cdots\ge|r_l|$.

Let $f_1\nek f_{2n}$ denote the dual basis of $V$.
\rem{cFL}
If the vector space $V$ is endowed with a complex structure $J:V\to V$
compatible with the metric (i.e., $J^*=-J$) and such that $F$ is a $(1,1)$ form
with respect to $J$, then the basis $f_1\nek f_{2n}$ can be chosen so that
$f_{j+n}=Jf_j, \ i=1\nek n$.
\erem

Let us define a complex structure $J:V\to V$ on $V$ by the condition
$f_{i+n}=Jf_i, \ i=1\nek n$. Then, the complexification of $V$ splits into the
sum of its holomorphic and anti-holomorphic parts
\[
        V\otimes\CC \ = \ V^{1,0}\oplus V^{0,1},
\]
on which $J$ acts by multiplication by $i$ and $-i$ respectively. The space
$V^{1,0}$ is spanned by the vectors $e_j=f_j-if_{j+n}$, and the space $V^{0,1}$
 is spanned by the vectors  $\oe_j=f_j+if_{j+n}$. Let  $e^1\nek e^n$
and $\oe^1\nek \oe^n$ be the corresponding dual base of $(V\ha)^*$ and
$(V\ah)^*$ respectively. Then \refe{iF} may be rewritten as
\[
        iF_x^\L \ = \  \frac{i}2\, \sum_{j=1}^n\, r_j\, e^j\wedge \oe^j.
\]
We will need the following simple
\lem{mur}
Let $\mu_1\nek\mu_l$ and $r_1\nek r_l$ be as above. Then $\mu_i=|r_i|$, for any
$i=1\nek l$. In particular,
\[
        \Tr^+\tJ \ = \ |r_1| +\cdots+|r_l|.
\]
\elem
\prf
Clearly, the vectors $e_1\nek e_n; \oe_1\nek\oe_n$ form a basis of eigenvectors
of $\tJ$ and
\begin{align}
        \tJ\, e_j \ &= \ ir_j\, e_j, \quad  \tJ\, \oe_j \ = \ -ir_j\, \oe_j \qquad
                        &\mbox{for} \quad j&=1\nek l, \notag\\
        \tJ\, e_j \ &= \ \tJ\, \oe_j \ = \ 0
                        \qquad &\mbox{for} \quad j&=l+1\nek n.\notag
\end{align}
Hence, all the nonzero eigenvalues of $\tJ$ are $\pm{i}|r_1|\nek\pm{i}|r_l|$.
\eprf

\ssec{spin1}{Spinors}
Set
\begin{equation}\label{E:spin1}
        S^+ \ = \ \bigoplus_{j \ \text{even}}\Lam^j(V\ah), \quad
         S^- \ = \ \bigoplus_{j \ \text{odd}}\Lam^j(V\ah).
\end{equation}
Define a graded action of the Clifford algebra $C(V)$ on the graded space
$S=S^+\oplus S^-$ as follows (cf. \refss{clstr}): \ if $v\in V$ decomposes as
 $v=v\ha+v\ah$  with  $v\ha\in V\ha$  and  $v\ah\in V\ah$, then its
Clifford action on $\alp\in E$ equals
\begin{equation}\label{E:clact}
        c(v)\alp \ = \ \sqrt{2}\, \left( v\ah\wedge\alp  \ - \ \iot(v\ha)\, \alp\right).
\end{equation}
Then (cf. \cite[\S3.2]{BeGeVe}) $S$ is the {\em spinor representation \/} of
$C(V)$, i.e., the complexification $C(V)\otimes\CC$ of $C(V)$ is isomorphic to
$\End(S)$. In particular, the Clifford module $E$ can be decomposed as
\[
        E \ = \ S\otimes W,
\]
where $W=\Hom_{C(V)}\, (S,E)$. The action of $C(V)$ on $E$ is equal to $a\mapsto
c(a)\otimes1$, where $c(a) \ (a\in C(V))$ denotes the action of $C(V)$
on $S$. The natural grading on $E$ is given by $E^\pm= S^\pm\otimes{W}$.

To prove \refp{cFL} it suffices now to study the action of $\bfc(F)$ on $S$. The
latter action is completely described be the following
\lem{ecF}
The vectors $\oe^{j_1}\wedge\dots\wedge\oe^{j_m}\in S$ form a basis of
eigenvectors of $\bfc(F)$ and
\[
        \bfc(F)\, \oe^{j_1}\wedge\dots\wedge\oe^{j_m} \ = \
                \big(\sum_{j'\not\in\{j_1\nek j_m\}} \, r_{j'} \ - \
                        \sum_{j''\in\{j_1\nek j_m\}} \, r_{j''}
                \big)\, \oe^{j_1}\wedge\dots\wedge\oe^{j_m}.
\]
\elem
\prf
The proof is an easy computation which is left to the reader.
\eprf

\ssec{prcFL}{Proof of \refp{cFL}}
Recall that the orientation of $V$ is fixed and that we have chosen the basis
$f_1\nek f_{2n}$ of $V$ which defines the same orientation. Suppose now that the
bilinear form $F$ is non-degenerate. Then $l=n$ in \refe{iF}. It is clear, that
the orientation defined by $iF$ coincides with the given orientation of $V$ if
and only if the number $q$ of negative numbers among $r_1\nek r_n$ is even.
Hence, by \refl{ecF}, the restriction of $\bfc(F)\in\End(S)$ on
$\Lam^j(V\ah)\subset S$ is greater than $-(\tau-2m)$ if the parity of $j$ and
$q$ are different. The \refp{cFL} follows now from \refe{spin1}.
\hfill$\square$

\ssec{prD-Dac}{Proof of \refp{D-Dac}}
Assume that at least $q$ of the numbers $r_1\nek r_l$ are negative and let
$m_q>0$ be the minimal positive number such that at least $q$ of these numbers
are not greater than $-m_q$. It follows from \refl{ecF}, that
\begin{equation}\label{E:>>>}
        \<\, c(F)\alp,\alp\, \> \ \ge \ -\, (\tau-2m_q)\, \|\alp\|^2,
\end{equation}
for any $j<q$ and any $\alp\in \Lam^j(V\ah)$.

Clearly, in the conditions of \refp{D-Dac}, the operator $A$ defined in \refe{A}
is $\G$-invariant, and, hence, bounded. \refp{D-Dac} follows now from \refe{>>>}
and the Lichnerowicz formula \refe{lichn}.
\hfill $\square$

\ssec{prgr}{Proof of \refp{gr}}
Let $\pi_i:L^2\A^{0,*}(M,\E\otimes\Lk)\to L^2\A^{0,i}(M,\E\otimes\Lk)$ denote
the projection and set 
\[
        \tilnk \ = \ \sum_i\, \pi_i\circ\nk\circ\pi_i.
\]
Denote $\tildk=(\tilnk)^*\tilnk$. 
Clearly, $\tildk$ preserves the $\ZZ$-grading on $L^2\A^{0,*}(M,\E\otimes\Lk)$.
It follows from the proof of Theorem~2.16 in
\cite{Dem85}, that there exists a sequence $\eps_1,\eps_2,\ldots$, convergent to
zero, such that
\[
        (1-\eps_k)\, \<\, \tildk\, \gam,\, \gam\, \> \ - \ \eps_kk\, \|\gam\|^2
        \ \le \
        \<\, \Del_k\, \gam,\,  \gam\, \>
        \ \le \
        (1+\eps_k)\, \<\, \tildk\, \gam,\, \gam\, \>
                \ + \ \eps_kk\, \|\gam\|^2,
\]
for any $\gam$ in the domain of $\Del_k$. Hence,
\begin{multline}\notag
        \|\<\, (\Del_k-\tildk)\, \gam,\,  \gam\, \>\| \\
        \ \le \
        \eps_k\, \<\, \tildk\, \gam,\, \gam\, \> \ + \ \eps_kk\, \|\gam\|^2
        \ \le \
        \eps_k \, \sum_{i}\, \<\, \Del_k\, \pi_i\gam,\,  \pi_i\gam\, \>
         \ + \ \eps_kk\, \|\gam\|^2.
\end{multline}

Suppose now that $\gam=\alp+\bet$, where
$\alp\in\Dom^i(D^2_k),\bet\in\Dom^j(D^2_k), \ i\not=j$. Then
\begin{multline}\notag
         \|\, 2\RE\<\, \Del_k\, \alp,\,  \bet\, \>\, \|
        \ = \
         \|\, 2\RE\<\, (\Del_k-\tildk)\, \alp,\,  \bet\, \>\, \| \\
        \ \le \
         \|\, \<\, (\Del_k-\tildk)\, \gam,\,  \gam\, \>\, \|
        \ + \
         \|\, \<\, (\Del_k-\tildk)\, \alp,\,  \alp\, \>\, \|
        \ + \
         \|\, \<\, (\Del_k-\tildk)\, \bet,\,  \bet\, \>\, \| \\
        \ \le \
         2\eps_k\, \<\, \Del_k\, \alp,\,  \alp\, \>
        \ + \
         2\eps_k\, \<\, \Del_k\, \bet,\,  \bet\, \>
        \ + \
         2\eps_kk\, \|\alp\|^2 \ + \ 2\eps_kk\, \|\bet\|^2.
\end{multline}
Similarly one obtains an estimate for the imaginary part of
$\<\Del_k\alp,\bet\>$. This leads to the following analogue of \refp{gr} for the
operator $\Del_k$:
\begin{equation}\label{E:albe}
         \|\, \<\, \Del_k\, \alp,\,  \bet\, \>\, \|
        \ \le  \
         2\eps_k\, \<\, \Del_k\, \alp,\,  \alp\, \>
        \ + \
         2\eps_k\, \<\, \Del_k\, \bet,\,  \bet\, \>
        \ + \
         2\eps_kk\, \|\alp\|^2 \ + \ 2\eps_kk\, \|\bet\|^2.
\end{equation}

We now apply the Lichnerowicz formula \refe{lechnW} to obtain \refp{gr} from
\refe{albe}. Note, first, that the operator
$A\in\End(\E)\subset\End(\E\otimes\Lk)$, defined in \refe{A}, is independent of
$k$ and bounded. Note, also, that, by \refl{ecF}, the operator $\bfc(F^\L)$
preserves the $\ZZ$-grading on $L^2\A^{0,*}(M,\Ek)$. Hence, it follows from
\refe{albe} and the Lichnerowicz formula \refe{lechnW} that
\begin{multline}\notag
          \|\, \<\, D_k^2\, \alp,\,  \bet\, \>\, \|
         \ \le \
          |\, \<\, \Del_k\, \alp,\,  \bet\, \>\, |
          \ + \
                   |\, \<\, A\, \alp,\,  \bet\, \>\, | \\
         \ \le \
         2\eps_k\, \<\, \Del_k\, \alp,\,  \alp\, \>
         \ + \
         2\eps_k\, \<\, \Del_k\, \bet,\,  \bet\, \>
         \ + \
         2\eps_kk\, \|\alp\|^2 \ + \ 2\eps_kk\, \|\bet\|^2
         \ + \
         \|A\|\, \|\alp\|\, \|\bet\| \\
        \ \le \
         2\eps_k\, \<\, D_k^2\, \alp,\,  \alp\, \>
         \ + \
         2\eps_k\, \<\, D_k^2\, \bet,\,  \bet\, \>
         \ + \
         2\eps_kk\, \big(1+\|\bfc(F^\L)\|+2\|A\|\big)\, \|\alp\|^2 \\
         \ + \
         2\eps_kk\, \big(1+\|\bfc(F^\L)\|+2\|A\|\big)\, \|\bet\|^2.
\end{multline}
Hence, \refp{gr} holds with $\del_k=(1+\|\bfc(F^\L)\|+2\|A\|\big)\eps_k$.
\hfill$\square$

\sec{lapl}{Estimate of the metric Laplacian}

In this section we prove \refp{lapl}. The proof consists of two steps. First, we
establish the following
\lem{compact}
Suppose that $F^\L$ does not vanish at any point $x\in M$. Fix a compact subset
$K\in{}M$ and a positive number $\eps>0$. There exists a constant $C_{K,\eps}$
such that, for any $k\in\ZZ$ and any smooth section $s{\in}\G(M,\Ek)$ supported
on $K$, one has
\eq{lapl-com}
    \<\Del_k\, s,s\> \ \ge \ k\, \<(\tau(x)-\eps) s,s\> \ - \ C_{K,\eps}\, \|s\|^2.
\end{equation}
\elem
Then we use the the fact $\Del_k$ is $\G$-invariant and that $M/\G$ is compact
to deduce \refp{lapl} from \refl{compact}. This is done in
Subsections~\ref{SS:IMS}, \ref{SS:prlapl}.

We now pass to the proof of \refl{compact}. The proof is almost verbatim
repetition of the proof of Proposition~4.4 in \cite{BrVanish} (see also
\cite[\S2]{GuiUr88}, \cite[\S2]{BoUr96}) and occupies
Subsections~\ref{SS:reduc}--\ref{SS:melin}.

\ssec{reduc}{Reduction to a scalar operator}
In this subsection we construct a space $\Z$ and an operator $\tilDel$ on the
space $L^2(\Z)$ of $\Z$, such that the operator $\Del_k$ is ``equivalent" to a
restriction of $\tilDel$ onto certain subspace of $L^2(\Z)$. This allow to
compare the operators $\Del_k$ for different values of $k$.

Let $\F$ be the principal $G$-bundle with a compact structure group $G$,
associated to the vector bundle $\E\to M$. Let $\Z$ be the principal
$(S^1\times{G})$-bundle over $M$, associated to the bundle $\E\otimes\L\to M$.
Then $\Z$ is a principle $S^1$-bundle over $\F$. We denote by $p:\Z\to \F$ the
projection.

The connection $\n^\L$ on $\L$ induces a connection on the bundle
$p:\Z\to\F$. Hence, any vector $X\in T\Z$ decomposes as a sum
\begin{equation}\label{E:dec}
        X \ = \ X{\hor} \ + \ X\vert,
\end{equation}
of its horizontal and vertical components.

Consider the {\em horizontal exterior derivative \/} $d\hor:
C^\infty(\Z)\to \A^1(\Z,\CC)$, defined by the formula
\[
        d\hor f(X) \ = \ df(X\hor), \qquad X\in T\Z.
\]

The connections on $\E$ and $\L$, the Riemannian metric on $M$, and
the Hermitian metrics on $\E, \L$ determine a natural Riemannian
metrics $g^\F$ and $g^\Z$ on $\F$ and $\Z$ respectively,
cf. \cite[Proof of Theorem~2.1]{BoUr96}. Let $(d\hor)^*$ denote the
adjoint of $d\hor$ with respect to the scalar products induced by this
metric. Let
\[
        \tD \ = \ (d\hor)^*d\hor:\, C^\infty(\Z) \to C^\infty(\Z)
\]
be the {\em horizontal Laplacian} for the bundle $p:\Z\to \F$.

Let $C^\infty(\Z)_k$ denote the space of smooth functions on $\Z$, which are
homogeneous of degree $k$ with respect to the natural fiberwise circle action on
the circle bundle $p:\Z\to \F$. It is shown in \cite[Proof of
Theorem~2.1]{BoUr96}, that to prove \refl{compact} it suffices to prove
\refe{lapl} for the restriction of $\tD$ to the space $C^\infty(\Z)_k$.

\ssec{symbol}{The symbol of $\tD$}
The decomposition \refe{dec} defines a splitting of the cotangent bundle $T^*\Z$
to $\Z$ into the horizontal and vertical subbundles. For any $\xi\in T^*\Z$, we
denote by $\xi\hor$ the horizontal component of $\xi$. Then, one easily checks
(cf. \cite[Proof of Theorem~2.1]{BoUr96}), that the principal symbol
$\sig_2(\tD)$ of $\tD$ may be written as
\begin{equation}\label{E:symb}
        \sig_2(\tD)(z,\xi) \ = \ g^\F(\xi\hor,\xi\hor).
\end{equation}
The subprincipal symbol of $\tD$ is equal to zero.

On the {\em character set \/} $\calC=\big\{(z,\xi)\in T^*\Z\backslash\{0\}:\,
\xi\hor=0\big\}$ the principal symbol $\sig_2(\tD)$ vanishes to second order.
Hence, at any point $(z,\xi)\in\calC$, we can define the {\em Hamiltonian map
\/} $F_{z,\xi}$ of $\sig_2(\tD)$, cf. \cite[\S21.5]{Horm3}. It is a
skew-symmetric endomorphism of the tangent space $T_{z,\xi}(T^*\Z)$. Set
\[
        \Tr^+ F_{z,\xi} \ = \ \nu_1 \ + \ \cdots \ + \ \nu_l,
\]
where $i\nu_1\nek i\nu_l$ are the nonzero eigenvalues of $F_{z,\xi}$
for which $\nu_i>0$.

Let $\rho:\Z\to M$ denote the projection. Then, cf. \cite[Proof of
Theorem~2.1]{BoUr96}
\footnote{The absolute value sign of $\xi\vert$ is erroneously missing
in \cite{BoUr96}.},
\begin{equation}\label{E:tr+F}
        \Tr^+ F_{z,\xi} \ = \ \tau(\rho(z))\, |\xi\vert|
\end{equation}
Here   $\xi\vert$ is the vertical component of $\xi\in T^*\Z$, and
$\tau$ denotes the function defined in \refe{tau}.

\ssec{melin}{Application of the Melin  inequality. Proof of \refl{compact}}
Let $D\vert$ denote the generator of the $S^1$ action on $\Z$. The symbol of
$D\vert$ is $\sig(D\vert)(z,\xi) =\xi\vert$. Fix $\eps>0$, and consider the
operator
\[
        A \ = \ \tD \ - \big(\tau(\rho(z))-\eps\big)\, D\vert:\,
                C^\infty(\Z) \ \to \  C^\infty(\Z).
\]
The principal symbol of $A$ is given by \refe{symb}, and the subprincipal symbol
\[
        \sig_1^s(A)(z,\xi) \ = \  -\big(\tau(\rho(z))-\eps\big)\, \xi\vert.
\]
It follows from \refe{tr+F}, that
\[
        \Tr^+ F_{z,\xi}+\sig_1^s(A)(z,\xi) \ \ge \
                \eps\, |\xi\vert| \ > \ 0.
\]
Hence, by the Melin inequality (\cite{Melin71}, \cite[Theorem~22.3.3]{Horm3}),
for any compact subset $K\subset \calZ$, there exists a constant $C_{K,\eps}$,
depending on $K$ and $\eps$, such that
\begin{equation}\label{E:Aff}
        \<\, A f,f\, \> \ \ge -\, C_{K,\eps}\, \|f\|^2,
                \quad \mbox{for any}\quad
                        f\in C^\infty(\Z), \ \supp f\subset K.
\end{equation}
Here $\|\cdot\|$ denotes the $L^2$ norm of the function $f\in{}C^\infty(\Z)$.

{}From \refe{Aff}, we obtain
\[
        \<\, \tD f,f\, \> \ \ge \
                \<\, (\tau(\rho(z))-\eps)D\vert f,f\, \> \
        - \  C_\eps\, \|f\|^2.
\]
Noting that if $f\in C^\infty(\Z)_k$, then $D\vert f=kf$, the proof of
\refl{compact} is complete.
\hfill $\square$

We pass to the proof of \refp{lapl}.
\ssec{IMS}{IMS localization formula}
If $f\in C^\infty(M)$ and $\gam\in \G$ we denote by $f^\gam$ the function
defined by the formula $f^\gam(x)=f(\gam^{-1}\cdot x)$. One easily sees (cf.
\cite[\S3]{Atiyah76}), that there exists a $C^\infty$ function $f:M\to [0,1]$
with compact support, such that
\begin{equation}\label{E:part}
        \sum_{\gam\in\G}\, \big(f^\gam(x)\big)^2 \ \equiv \ 1.
\end{equation}
We need the following version of the IMS localization formula (cf. \cite{CFKS},
\cite[Lemma~4.10]{BrFar1})
\lem{IMS}
The following identity holds
\begin{equation}\label{E:IMS}
        \Del_k \ = \ \sum_{\gam\in\G}\, f^\gam\Del_k f^\gam  \ + \
        \frac12\, \sum_{\gam\in\G}\, [f^\gam,[f^\gam,\Del_k]].
\end{equation}
\elem
\prf
Using \refe{part}, we can write
\[
        \Del_k \ = \ \sum_{\gam\in\G}\, (f^\gam)^2\Del_k \ = \
          \sum_{\gam\in\G}\, \Big(\,
            f^\gam\Del_k f^\gam + f^\gam[f^\gam,\Del_k] \, \Big).
\]
Similarly,
\[
        \Del_k \ = \ \sum_{\gam\in\G}\, \Del_k(f^\gam)^2 \ = \
          \sum_{\gam\in\G}\, \Big(\,
            f^\gam\Del_k f^\gam - [f^\gam,\Del_k]f^\gam \, \Big).
\]
Summing these two identities and dividing by 2 we come to \refe{IMS}.
\eprf

\ssec{prlapl}{Proof of \refp{lapl}}
To prove \refp{lapl} it remains to estimate all the terms in the equality
\refe{IMS}.

It follows from \refl{compact}, that for any $\eps>0$ there exists a constant
$C_\eps>0$ such that
\[
        \<f\Del_k f\, s, s\> \ = \ \<\Del_k\, f s,f s\> \ \ge \
                k\, \<(\tau(x)-\eps) f^2s,s\> \ - \ C_{\eps}\, f^2\|s\|^2,
                \qquad s\in\G(M,\Ek).
\]
Since the operator $\Del_k$ is $\G$-invariant, this inequality remains true if
we replace everywhere $f$ with $f^\gam$. Hence, in view of \refe{part},
\begin{equation}\label{E:est1}
        \sum_{\gam\in\G}\,  \<f^\gam\Del_k f^\gam\, s, s\> \ \ge \
                k\, \<(\tau(x)-\eps) s,s\> \ - \ C_{\eps}\, \|s\|^2.
\end{equation}

We now study the second summand in \refe{IMS}. Since $\Del_k$ is a second order
differential operator, the double-commutant $[f,[f,\Del_k]]$ is an operator of
multiplication by a function. Let us denote this function by $\Phi$. Then, for
any $\gam\in\G$, the operator $[f^\gam,[f^\gam,\Del_k]]$ acts by multiplication
by $\Phi^\gam$. Since the support of $f$ is compact so is the support of $\Phi$.
It follows that
\begin{equation}\label{E:est2}
        \sum_{\gam\in\G}\, [f^\gam,[f^\gam,\Del_k]] \ = \
            \sum_{\gam\in\G}\, \Phi^\gam
\end{equation}
is a smooth $\G$-invariant function. Hence, it is bounded. Set
\[
        C \ = \ \max_{x\in M}\, \big|\, \sum_{\gam\in\G}\, \Phi^\gam\, \big|.
\]
From \refe{est2} we obtain
\begin{equation}\label{E:est3}
        \big\|\, \sum_{\gam\in\G}\, [f^\gam,[f^\gam,\Del_k]]\, \big\| \ \le \ C.
\end{equation}
Combining \refe{IMS}, \refe{est1} and \refe{est3} we obtain
\[
        \<\Del_k \, s, s\> \ = \ \<\Del_k\, s, s\> \ \ge \
                k\, \<(\tau(x)-\eps) f^2s,s\> \ - \ (C_{\eps}+\frac12C)\,
                f^2\|s\|^2,
\]
for any $s\in \G(M,\Ek)$. \hfill$\square$

\providecommand{\bysame}{\leavevmode\hbox to3em{\hrulefill}\thinspace}

\end{document}